\magnification = 1100
\documentstyle{amsppt} 
\overfullrule0pt
\vsize=53 true pc
\TagsOnLeft  
\def \n {\noindent}
\def \s {\smallskip}
\def \m {\medskip}
\def \b {\bigskip}
\def \a {\alpha}
\def \ad {{\operatorname {ad}}}
 
\def \be {\beta} 

\def \cop {{\operatorname {coop}}}

\def \e {\epsilon} 
\def \End {{\operatorname {End}}}
  
\def\fsl{\frak s\frak l_2} 
\def\fsn{\frak s\frak l_n}
\def\fgl{\frak g\frak l_n}

\def \ga {\gamma}  
\def \g {{\frak g}}

\def \half {\frac{1}{2}}

\def \ii {{\operatorname {(i)}}}

\def \l {\lambda}

\def \la {\langle}

\def \op {{\operatorname {op}}} 
\def \ot {\otimes}

\def \ra {\rangle}

\def \ve {\check {e}}
\def \vf {\check {f}}
\def \vo {\check {\omega}}
\def \va {\check {a}}
\def \vb {\check {b}} 
\def \w {\omega}

\def \K {\Bbb K} 
\def \Z {\Bbb Z}
\font\bigf=cmbx10 scaled \magstep1

\document
\baselineskip=12pt 
\topmatter
\title {\bigf Two-Parameter Quantum Groups \\
 and Drinfel'd Doubles} \endtitle
\author Georgia Benkart \footnote {The authors gratefully acknowledge 
support from  National Science
Foundation Grant \#{}DMS--9970119, National Security Agency Grant
\#{}MDA904-01-1-0067, and the hospitality of
the Mathematical Sciences Research Institute, Berkeley. \quad\quad 
\quad\quad \quad\quad \quad\quad\quad\quad\quad\quad\quad} \\ 
Sarah Witherspoon$^1$
\endauthor 
\leftheadtext{ GEORGIA BENKART, \; SARAH WITHERSPOON} \rightheadtext 
{TWO-PARAMETER QUANTUM GROUPS}
\date August 3, 2001 
\enddate \subjclass\nofrills \n 2000 {\it Mathematics 
Subject Classification.\/} Primary 17B37, 16W30, 16W35, 81R50 
\endsubjclass \abstract We investigate two-parameter quantum groups 
corresponding to the general linear and special linear Lie algebras 
$\fgl$ and $\fsn$.  We show that these quantum groups can be realized 
as Drinfel'd doubles of certain Hopf subalgebras with respect to Hopf 
pairings.  Using the Hopf pairing, we construct a corresponding 
$R$-matrix and a quantum Casimir element.  We discuss isomorphisms 
among these quantum groups and connections with multiparameter quantum 
groups.  \endabstract \endtopmatter \head Introduction \endhead \m In 
this work we study two two-parameter quantum groups $\widetilde U = 
U_{r,s}(\fgl)$ and $U = U_{r,s}(\fsn)$ corresponding to the Lie 
algebras $\fgl$ and $\fsn$.  Our Hopf algebra $\widetilde U$ is 
isomorphic as an algebra to Takeuchi's $U_{r,s^{-1}}$ (see [T]), but 
as a Hopf algebra, it has the opposite coproduct.  (A different 
presentation of $U_{r,s^{-1}}$ was obtained by Kulish [K] (see also  
[Ji]).) As an algebra, $\widetilde U$ has generators $e_j, \ f_j, (1 
\leq j < n)$, and $\ a_i^{\pm 1}, \ b_i^{\pm 1}$ \ \ ($1 \leq i \leq 
n$), and defining relations given in (R1)-(R7) below.  The elements 
$e_j, \ f_j$, \ $\w_j^{\pm 1}, (\w_j')^{\pm 1}$ ($1 \leq j < n)$, 
where $\w_j = a_jb_{j+1}$ and $\w_j' = a_{j+1}b_j$, generate the 
subalgebra $U = U_{r,s}(\fsn)$.  \m We show that both $\widetilde{U}$ 
and $U$ may be realized as Drinfel'd doubles of certain Hopf 
subalgebras with respect to suitable Hopf pairings.  
Using the Hopf pairing, we construct an 
$R$-matrix for $\widetilde{U}$ (which also works for $U$).  For
$\widetilde{U}$-modules $M$ and $M'$ in
category $\Cal O$ (defined in Section 4), there is an isomorphism $R_{M',M}: 
M'\otimes M\rightarrow M\otimes M'$.  Moreover, the $R$-matrix 
satisfies the quantum Yang-Baxter equation and the hexagon identities.  
In [BW2], the $R$-matrix will be used to establish an analogue of 
Schur-Weyl duality in this setting: \ $\widetilde{U}$ has a natural 
$n$-dimensional module $V$, and the centralizer algebra 
$\End_{\widetilde{U}}(V^{\ot k})$ is generated by a certain Hecke 
algebra $H_k(r,s)$.  We construct a quantum Casimir element, which 
will play an essential role in [BW2] in proving that 
finite-dimensional modules in category $\Cal O$ are completely 
reducible.  \m Jing's work [Ji], which treats the special case of 
$\frak {gl}_2$, adopts exactly the opposite approach to the one of 
this paper $-$ it derives an analogue of the algebra $\widetilde U$ 
from one particular solution $R$ of the quantum Yang-Baxter equation.  
Similarly, Chin and Musson [ChM] and Dobrev and Parashar [DP] study 
multiparameter quantum universal enveloping algebras defined as duals 
of quantum function algebras arising from $R$-matrices.  In Section 6, 
we relate the two-parameter quantum groups considered here with 
certain special cases of these multiparameter quantum groups.  
Moreover, we determine conditions for isomorphisms among the 
two-parameter quantum groups.  In particular, the standard 
one-parameter quantum group $U_q(\fsl)$ of [Ja] is isomorphic to a 
quotient of $U_{r,s}(\fsl)$ by the ideal generated by $\omega_1' - 
\omega_1^{-1}$ whenever $q$ is a square root of $rs^{-1}$.  However, 
for $n \geq 3$, no such isomorphism exists (see Proposition 6.1).  \m 
Our motivation to study these two-parameter quantum groups came from 
our work [BW1] on down-up algebras.  Down-up algebras were introduced 
in [BR] as a generalization of the algebra generated by the down and 
up operators on posets.  They are unital associative algebras 
$A(\a,\be,\ga)$ over a field $\K$ having generators $d,u$ which 
satisfy the defining relations

$$\aligned
d^2u = \alpha dud + \beta u d^2 + \gamma d \\
d u^2 = \alpha udu + \beta u^2 d + \gamma u, \\
\endaligned$$

\n where $\a,\be,\gamma$ are fixed but arbitrary scalars in $\K$. 
If $\gamma \neq 0$, then the down-up algebra
$A(\alpha,\beta,\gamma)$ is isomorphic to
$A(\alpha,\beta,1)$.   Thus, there are basically
two different cases: $\gamma = 0$ and $\gamma = 1$.
Examples of down-up algebras include the universal
enveloping algebras of $\fsl$, of the Heisenberg Lie
algebra, and of the Lie superalgebra $\frak{osp}(1,2)$, which are
$A(2,-1,1)$, $A(2,-1,0)$, and  $A(0,-1,1)$, respectively, and many
of Witten's deformations of $U(\fsl)$ (see [B]).
Down-up algebras exhibit
many striking features including a Poincar\'e-Birkhoff-Witt 
type basis and a well-behaved representation theory ([BR], [KMP],
[CaM], [Jor], [KK1], [KK2], [Ku], [BL, Sec. 4]).   They are Noetherian domains 
whenever $\be \neq 0$. 
 \m
Essential to the structure of $A(\a,\be,\ga)$ are the roots 
of the equation 

$$0 = t^2 - \a t -\be = (t-r)(t-s),$$

\n Thus, $\a = r+s$ and $\be = -rs$.  
When $rs \neq 0$ and $\ga = 0$, the down-up algebra $A(\a,\be,0) =A(r+s,-rs,0)$ can be
extended by automorphisms to give a  Hopf algebra $B(r+s,-rs,0)$ (see [BW1]).  This Hopf
algebra is isomorphic to a subalgebra of $U_{r,s}(\frak{sl}_3)$ when
$r$ and $s$ are not roots of unity and to a quotient of a subalgebra 
when they are.
It seemed natural to expect that there is a Drinfel'd double 
(quantum double) of the subalgebra, which
yields a quantum group that depends on the two parameters $r$ and $s$.   In fact, that
quantum group is $U_{r,s}(\frak{sl}_3)$. That result is a very special case of our theorem
showing that $U$ and $\widetilde U$ are Drinfel'd doubles.    
\m
Throughout we will be working over a field $\K$, which is required to be
algebraically closed from Section 3 to the end of the
paper.
\b \m
 \head \S 1.  Preliminaries \endhead \m Assume $\Phi$ is a finite 
root system of type A$_{n-1}$ with $\Pi$ a base of simple roots.  We 
regard $\Phi$ as a subset of a Euclidean space $E = \Bbb R^n$ with an 
inner product $\la \,,\,\ra$.  We let $\e_1, \dots, \e_n$ denote an 
orthonormal basis of $E$, and suppose $\Pi = \{\a_j = \e_{j}- \e_{j+1} 
\mid j = 1, \dots, n-1\}$ and $\Phi = \{\e_i -\e_j \mid 1 \leq i \neq 
j \leq n\}$.  \m Fix nonzero elements $r,s$ in a field $\K$.  Here we 
assume $r \neq s$.  \m Let $\widetilde U = U_{r,s}(\fgl)$ be the 
unital associative algebra over $\K$ generated by elements $e_j, \ 
f_j, (1 \leq j < n)$, and $\ a_i^{\pm 1}, \ b_i^{\pm 1}$ \ \ ($1 \leq 
i \leq n$), which satisfy the following relations.  \m \roster 
\item"{(R1)}" The $a_i^{\pm 1}, \ b_j^{\pm 1}$ all commute with one 
another and \ $a_ia_i^{-1}= b_j b_j^{-1}=1,$ \m \item"{(R2)}" $ a_ie_j 
= r^{\la\epsilon_i,\a_j\ra}e_j a_i$ \ \ and \ \ $a_if_j 
=r^{-\la\epsilon_i,\a_j\ra} f_ja_i,$ \m \item"{(R3)}" $b_ie_j = 
s^{\la\epsilon_i,\a_j\ra}e_j b_i$ \ \ and \ \ $b_if_j =s^{-\la 
\epsilon_i,\a_j\ra} f_jb_i,$ \m \item"{(R4)}" 
$\displaystyle{[e_i,f_j]=\frac{\delta_{i,j}}{r-s}(a_ib_{i+1}-a_{i+1}b_i),}$ 
\m \item"{(R5)}" $[e_i,e_j]=[f_i,f_j]=0 \ \ \text{ if }\ \ |i-j|>1, $ 
\m \item"{(R6)}" $e_i^2e_{i+1}-(r+s)e_ie_{i+1}e_i+rse_{i+1}e_i^2 = 0,$ 
\smallskip \hskip -.2 truein $e_i e^2_{i+1} -(r+s)e_{i+1}e_ie_{i+1} 
+rs e^2_{i+1}e_i = 0,$ \m \item"{(R7)}" 
$f_i^2f_{i+1}-(r^{-1}+s^{-1})f_if_{i+1}f_i +r^{-1}s^{-1}f_{i+1}f_i^2 = 
0,$ \smallskip \hskip -.2 truein $f_i f^2_{i+1} 
-(r^{-1}+s^{-1})f_{i+1}f_if_{i+1}+r^{-1}s^{-1} f^2_{i+1} f_i=0.$ 
\endroster \m The relations in (R6) are just the two defining 
relations of the down-up algebra $A(r+s,-rs,0)$, while those in (R7) 
are the defining relations of $A(r^{-1}+s^{-1},-r^{-1}s^{-1},0)$.  
In fact these two down-up algebras are isomorphic via the map that 
takes $d$ to $u'$ and $u$ to $d'$ (assuming $d',u'$ are the generators 
of the latter) (see [BR]).

\m
We will be interested in the subalgebra $U = U_{r,s}(\fsn)$ of $\widetilde U
= U_{r,s}(\fgl)$ generated by the elements $e_j,f_j$, $\w_j$, and
$\w_j'$ \ ($1 \leq j < n)$, where

$$\omega_j = a_j b_{j+1}  \ \ \text{and} \ \ \omega_j' = a_{j+1}b_j. \tag 1.1$$
\n 
These elements satisfy (R5)-(R7) along with the following relations:
\m
\roster
\item"{(R1')}"  The $\w_i^{\pm 1}, \ \w_j^{\pm 1}$ all commute with one 
another and $\w_i \w_i^{-1}= \w_j'(\w_j')^{-1}=1,$
\m
\item"{(R2')}"  $ \w_i e_j = r^{\la \e_i,\a_j\ra }s^{\la \e_{i+1},\a_j\ra}e_j \w_i$ \ \ and
\ \
$\w_if_j = r^{-\la\e_i,\a_j\ra}s^{-\la\e_{i+1},\a_j\ra} f_j\w_i,$ 
\m
\item"{(R3')}"
$\w_i'e_j = r^{\la\epsilon_{i+1},\a_j\ra}s^{\la\e_{i},\a_j\ra}e_j \w_i'$  \ \ and
\ \
$\w_i'f_j = r^{-\la\epsilon_{i+1},\a_j\ra}s^{-\la\e_{i},\a_j\ra}f_j \w_i'$,  
\m 
\item"{(R4')}"
$\displaystyle{[e_i,f_j]=\frac{\delta_{i,j}}{r-s}(\w_i-\w_i').}$
\endroster
\m 
When $r = q$ and $s = q^{-1}$, the algebra $U_{r,s}(\fgl)$ modulo
the ideal generated by the elements $b_i-a_i^{-1}$, $1 \leq i \leq n$,  is just the
quantum general linear group $U_q(\fgl)$,  and $U_{r,s}(\fsn)$ modulo
the ideal generated by the elements $\w_j'-\w_j^{-1}$, $1 \leq j < n$, is $U_q(\fsn)$.
\m
Let $Q = \Z \Phi$ denote the root lattice and set $Q^+ = \sum_{i = 1}^{n-1} \Z_{\geq 0}
\a_i$.  Then for  any $\zeta = \sum_{i=1}^{n-1}\zeta_i \a_i \in Q$, we  adopt the
shorthand 

$$
\w_\zeta   = \w_1^{\zeta_1} \cdots \w_{n-1}^{\zeta_{n-1}}, \ \ \ \ \
\w_\zeta' = (\w_1')^{\zeta_1} \cdots (\w_{n-1}')^{\zeta_{n-1}}
\tag 1.2$$

\n The following lemma is straightforward to check.
\b
\proclaim{Lemma 1.3} Suppose that $\zeta = \sum_{i=1}^{n-1}\zeta_i \a_i \in Q$. Then 
$$
\aligned w_\zeta e_i & = r^{-\la\e_{i+1},\zeta\ra}s^{-\la\e_i,\zeta\ra}e_i w_\zeta 
\qquad w_\zeta f_i = r^{\la\e_{i+1},\zeta\ra}s^{\la\e_i,\zeta\ra}f_i w_\zeta \\
w_\zeta' e_i & = r^{-\la\e_{i},\zeta\ra}s^{-\la\e_{i+1},\zeta\ra}e_i w_\zeta' 
\qquad w_\zeta' f_i = r^{\la\e_{i},\zeta\ra}s^{\la\e_{i+1},\zeta\ra}f_i w_\zeta'. 
\endaligned$$ \endproclaim 
\b
The algebras $\widetilde U$ and $U$ are Hopf algebras, 
where the $a_i^{\pm 1}, b_i^{\pm 1}$ are group-like elements, and the remaining
coproducts are determined by  

$$\Delta(e_i)=e_i\otimes 1 + \omega_i\otimes e_i, \quad  \ 
\Delta(f_i)=1\otimes f_i + f_i\otimes \omega_i'. \tag 1.4$$

\n This forces the counit and antipode maps to be

$$\aligned
& \varepsilon (a_i) = \varepsilon (b_i) = 1, \ \ \ \ S(a_i) = a_i^{-1},
\ \ \ S(b_i) = b_i^{-1} \\ 
& \varepsilon(e_i)=\varepsilon(f_i)=0, \ \ \
\, S(e_i)=-\omega_i^{-1}e_i, \ \ \ S(f_i)= -f_i(\omega_i')^{-1}.
\endaligned \tag 1.5$$

\b \m \head \S 2.  Drinfel'd doubles \endhead \m A {\it Hopf pairing} 
of two Hopf algebras $H$ and $H'$ is a bilinear form on $H'\times H$ 
satisfying the following properties (see [Jo, 3.2.1]): \m \roster 
\item"{(2.1)(i)}" $(1,h)=\varepsilon_H(h), \ \ \ (h',1) = 
\varepsilon_{H'}(h')$ \m \item"{(ii)}" 
$(h',hk)=(\Delta_{H'}(h'),h\otimes k) = \sum (h_{(1)}',h)(h_{(2)}',k)$ 
\m \item"{(iii)}" $(h'k',h)=(h'\otimes k',\Delta_H(h)) = \sum 
(h',h_{(1)}) (k', h_{(2)})$ \endroster \m \n for all $h,k\in H$ and 
$h',k'\in H'$, where $\varepsilon_H$ and $\varepsilon _{H'}$ denote 
the counits of $H$ and $H'$, respectively, and $\Delta_H$ and 
$\Delta_{H'}$ are their coproducts.  It is a consequence of the 
defining properties that a Hopf pairing satisfies
$$(S_{H'}(h'),h)=(h',S_{H}(h))$$
for all $h\in H$ and $h'\in H'$, where $S_H$ and $S_{H'}$ denote the respective antipodes
of $H$ and $H'$.
\m
Assume $\widetilde{B}$ is the Hopf subalgebra 
of $\widetilde{U}$ generated by $e_j, \ \omega_j^{\pm 1} \
(1\leq j< n)$, and $a_n^{\pm 1}$.  
Let $(\widetilde{B}')^{\cop}$ be the Hopf algebra having the
{\it opposite coproduct} to the Hopf subalgebra of $\widetilde{U}$ generated by
$f_j,\ (\omega_j')^{\pm 1} \ (1\leq j <n)$, and $b_n^{\pm 1}$.  
Similarly $B$ is the Hopf subalgebra of $U$ generated by $e_j, \omega_j^{\pm 1}
 \ (1\leq j<n)$, and $(B')^{\cop}$ is generated by $f_j,(\omega_j')^{\pm 1} 
\ (1\leq j <n)$.
\b
\proclaim{Lemma 2.2}
There are Hopf pairings of $\widetilde{B}$ and $\widetilde{B}'$,
respectively of $B$ and $B'$.
\endproclaim 
\m
\demo{Proof} We begin by defining a bilinear form for $\widetilde{B}'\times 
\widetilde{B}$ first on the generators:

$$\aligned
(f_i,e_j) &= \frac{\delta_{i,j}}{s-r},\\
(\omega_i',\omega_j) &= r^{\la\e_j,\a_i\ra}
s^{\la\epsilon_{j+1},\a_i\ra} = r^{-\la\e_{i+1},\a_j\ra}s^{-\la\e_i,\a_j\ra},\\
(b_n,a_n)&=1, \quad (b_n,\omega_j)=s^{-\la\epsilon_n,\alpha_j\ra}, \quad (\omega_i', a_n)
=r^{\la\epsilon_n,\alpha_i\ra},
\endaligned \tag 2.3$$

\n for $(1\leq i,j<n)$.
If $\omega_i'$ is replaced by $(\omega_i ')^{-1}$ in the second or third line
of (2.3), we replace the image under the bilinear form
by its inverse, and similarly for $\omega_j$ and $\omega_j^{-1}$, $a_n$ and
$a_n^{-1}$, $b_n$ and $b_n^{-1}$.
On all other pairs of generators the form is 0.
\m
In the second line of (2.3) we have applied the identity
$$\la \e_j,\a_i\ra = -\la\e_{i+1},\a_j\ra, \tag 2.4$$

\n which is quite useful in subsequent calculations. 
\m
The pairings in (2.3) may be {\it extended} to a bilinear form on 
$\widetilde{B}'\times \widetilde{B}$ by requiring that (2.1)(i)--(iii) hold.
We need only verify that the relations in $\widetilde{B}$ and $\widetilde{B}'$
are preserved, ensuring that the bilinear form is well-defined. 
It will then be a Hopf pairing by definition.  
Restricting the form to $B'\times B$ gives
the desired Hopf pairing of $B$ and $B'$.
\m
It is straightforward to check that
the bilinear form preserves all the relations among the $\omega_i^{\pm 1}, \
a_n^{\pm 1}$ in $\widetilde B$ and the $(\omega_j')^{\pm 1}, \ b_n^{\pm 1}$ in
 $\widetilde{B}'$.
We will verify that the form on $\widetilde{B}'\times \widetilde{B}$
preserves one of the remaining relations in $\widetilde{B}$, and leave the
other verifications to the reader. 
For each $i$, $1\leq i<n$, consider

$$(X, e_i^2e_{i+1}-(r+s)e_ie_{i+1}e_i+rse_{i+1}e_i^2),$$ 
where $X$ is any word in the generators
of $\widetilde{B}'$.  By definition, this is equal to

$$(\Delta^2(X), e_i\otimes e_i\otimes e_{i+1} -(r+s)
e_i\otimes e_{i+1}\otimes e_i +rse_{i+1}\otimes e_i\otimes e_i).\tag 2.5$$
In order for any one of these terms to be nonzero, $X$ must involve exactly
two $f_i$ factors, one $f_{i+1}$ factor, and arbitrarily many 
$(\omega_j')^{\pm 1}$ and $b_n^{\pm 1}$ factors $(1\leq j<n)$.  
First assume that 
$X=f_i^2f_{i+1}$. Then $\Delta^2(X)$ is equal to
$$(\omega_i'\otimes \omega_i'\otimes f_i + \omega_i'\otimes f_i\otimes 1 + 
f_i\otimes 1\otimes 1)^2(\omega_{i+1}'\otimes\omega_{i+1}'\otimes f_{i+1}
+ \omega_{i+1}'\otimes f_{i+1}\otimes 1 + f_{i+1}\otimes 1\otimes 1).$$
The relevant terms of $\Delta^2(X)$ are

$$\aligned
& f_i\omega_i'\omega_{i+1}'\otimes f_i\omega_{i+1}'\otimes f_{i+1}
+\omega_i'f_i\omega_{i+1}'\otimes f_i\omega_{i+1}'\otimes f_{i+1}\\
& \ \ \ \ \ +f_i\omega_i'\omega_{i+1}'\otimes \omega_i'f_{i+1}\otimes f_i
+\omega_i'f_i\omega_{i+1}'\otimes \omega_i'f_{i+1}\otimes f_i\\
& \ \ \ \ \ \ \ \ +(\omega_i')^2f_{i+1}\otimes f_i\omega_i'\otimes f_i
+(\omega_i')^2f_{i+1}\otimes \omega_i'f_i\otimes f_i.
\endaligned$$
Therefore (2.5) becomes

$$\aligned
 & (f_i\omega_i'\omega_{i+1}',e_i)(f_i\omega_{i+1}',e_i)(f_{i+1},e_{i+1})
+(\omega_i'f_i\omega_{i+1}',e_i)(f_i\omega_{i+1}',e_i)(f_{i+1},e_{i+1})\\
& \ \  -(r+s)(f_i\omega_i'\omega_{i+1}',e_i)(\omega_i'f_{i+1},e_{i+1})
(f_i,e_i) -(r+s)(\omega_i'f_i\omega_{i+1}',e_i)(\omega_i'f_{i+1},e_{i+1})
(f_i,e_i)\\
& \qquad \ \  +rs((\omega_i')^2f_{i+1},e_{i+1})(f_i\omega_i',e_i)(f_i,e_i)
+rs((\omega_i')^2f_{i+1},e_{i+1})(\omega_i'f_i,e_i)(f_i,e_i)\\
&\qquad = \frac{1}{(s-r)^3}( 1+(\omega_i',\omega_i) -(r+s)(\omega_i',
\omega_{i+1})-(r+s)(\omega_i',\omega_i)(\omega_i',\omega_{i+1})\\
& \qquad \qquad \quad  +rs(\omega_i',
\omega_{i+1})^2 +rs(\omega_i',\omega_{i+1})^2(\omega_i',\omega_i))\\
& \qquad = \frac{1}{(s-r)^3}(1+rs^{-1}-(r+s)r^{-1}-(r+s)s^{-1}+r^{-1}s +1) \ = \
0.
\endaligned$$

\n If $X=f_if_{i+1}f_i$ or $X=f_{i+1}f_i^2$, 
then similar calculations show that (2.5) is equal to 0.
Finally if $X$ is {\it any} word involving exactly two $f_i$ factors, 
one $f_{i+1}$
factor, and arbitrarily many factors of $(\omega_j')^{\pm 1} \
(1\leq j<n)$ and $b_n^{\pm 1}$, then (2.5) will just be a scalar multiple of 
one of the quantities we have already calculated, and therefore will equal 0.
(For example, if $X=f_i^2\omega_j'f_{i+1}$, then (2.5) will be 
$(\omega_j',\omega_{i+1})$ times the corresponding quantity for $X=
f_i^2f_{i+1}$.)

Analogous calculations show that the relations in $\widetilde{B}'$ 
are preserved.
\qed \enddemo
\b
As there is a Hopf pairing between $\widetilde{B}$ and $\widetilde{B}'$, 
there is a skew-Hopf pairing between $\widetilde{B}$ and 
$(\widetilde{B}')^{\cop}$, where the latter is $\widetilde{B}'$ as an algebra,
but with the opposite coproduct. Therefore, we may form the Drinfel'd double
$D(\widetilde{B},(\widetilde{B}')^{\cop})$ as in [Jo, 3.2].  This is a Hopf 
algebra whose underlying coalgebra is $\widetilde{B}\otimes (\widetilde{B}')
^{\cop}$ (that is, $\widetilde{B}\otimes (\widetilde{B}')^{\cop}$ as a
vector space with the tensor product coalgebra structure).  The algebra
structure is given as follows:  $\widetilde{B}$ and $\widetilde{B}'$ are 
identified as algebras with $\widetilde{B}\otimes 1$ and 
$1\otimes \widetilde{B}'$ respectively in $D(\widetilde{B},
(\widetilde{B}')^{\cop})$.  Letting $a\in \widetilde{B}$ and $b\in 
\widetilde{B}'$, we have $(a\otimes 1)(1\otimes b)= a\otimes b$ and
$$ (1\otimes b)(a\otimes 1) = \sum(S^{\cop}(b_{(1)}),a_{(1)})(b_{(3)},a_{(3)})
a_{(2)}\otimes b_{(2)},$$ 

\n where $S^{\cop}$ denotes the antipode for
$(\widetilde{B}')^{\cop}$.   (This expression looks different from  
[Jo, Lemma 3.2.2 (iii)] as we have written our bilinear form in the reverse 
order.)
A similar construction applies to $B$ and $B'$.
\b
\proclaim{Theorem 2.7}
$D(\widetilde{B},(\widetilde{B}')^{\cop})$ is isomorphic to $\widetilde{U}$,
and $D(B,(B')^{\cop})$ is isomorphic to $U$.
\endproclaim
\m
 \demo{Proof} We will prove the first statement, and the second will 
follow by restricting to fewer generators.  We will denote the image 
$e_i \ot 1$ of $e_i$ in $D(\widetilde{B},(\widetilde{B}')^{\cop})$ by 
$\ve_i$, and similarly for $\omega_i, a_n, f_i, \omega_i', b_n$.  
Define a map 
$\varphi:D(\widetilde{B},(\widetilde{B}')^{\cop})\rightarrow 
\widetilde{U}$ by $$\gathered \varphi(\ve_i)=e_i, \quad \quad 
\varphi(\vf_i)=f_i,\\
\varphi(\vo_i^{\pm 1})=\omega_i^{\pm 1},  \ \ 
\varphi((\vo_i')^{\pm 1})=   (\omega_i')^{\pm 1}, \   \ 
\varphi(\va_n^{\pm 1})=a_n^{\pm 1}, \ \ 
\varphi(\vb_n^{\pm 1})=b_n^{\pm 1}.
\endgathered$$
First notice that by definition, $\varphi$ preserves the coalgebra structures, 
the relations in $\widetilde{B}$, and the relations in $\widetilde{B}'$.
Next we will verify that the mixed relations in $D(\widetilde{B},
(\widetilde{B}')^{\cop})$ correspond to those in $\widetilde{U}$. 

To calculate $\vf_j \ve_i $, we use
$$\aligned
& \Delta^2(e_i)  = e_i\otimes 1\otimes 1 + \omega_i\otimes e_i\otimes 1
+\omega_i\otimes\omega_i\otimes e_i,\\
& (\Delta^{\cop})^2(f_j) = 1\otimes 1\otimes f_j + 1\otimes f_j\otimes\omega_j'
+f_j\otimes \omega_j'\otimes \omega_j',\\
& \text{and } \ \ S^{\cop}(f_j)  = -f_j(\omega_j')^{-1},
\endaligned$$
so that 
$$\aligned
\vf_j \ve_i  &= (-f_j(\omega_j')^{-1},e_i)(\omega_j',1)\vo_j'  
+ (1,\omega_i)
(\omega_j',1)\ve_i\vf_j  +(1,\omega_i)(f_j,e_i)\vo_i \\
&= -\frac{\delta_{i,j}}{s-r}\vo_j'  +\ve_i\vf_j  
+\frac{\delta_{i,j}}{s-r}\vo_i .
\endaligned$$
That is, $[\ve_i,\vf_j] =  \delta_{i,j}(s-r)^{-1}(\vo_i 
-\vo_i')$.  
Applying $\varphi$ gives the desired relation in $\widetilde{U}$.

We leave verification of the remaining relations to the reader.
As $\widetilde{U}$ is generated by $e_i, f_i, \omega_i^{\pm 1}, (\omega_i')^{\pm 1}
\ (1\leq i < n)$, $a_n$ and $b_n$, the map $\varphi$ is surjective, and there is an
obvious inverse map.
\qed\enddemo
\b 
\m  
\head {\S 3. Weight modules} \endhead 
\m
Let $\Lambda= \Z \e_1 \oplus \cdots \oplus \Z \e_n$, which is
the weight
lattice of $\fgl$.    Corresponding to any  $\l \in \Lambda$
is an algebra homomorphism $\hat \l$ from the subalgebra
$\widetilde U^{0}$ of $\widetilde U$ generated by the elements $a_i^{\pm
1}, b_i^{\pm 1}$ ($1 \leq i \leq n$)  
to $\K$ given by 

$$ \hat \l(a_i) =  r^{\la\e_i,\lambda\ra} \qquad \text{and} \qquad 
\hat \l(b_i) = s^{\la\e_{i},\lambda\ra}. \tag 3.1$$ 

\n The restriction $\hat \l: U^{0} \rightarrow \K$ to the subalgebra
$U^0$ of
$U$ generated by $\w_j^{\pm 1}, (\w_j')^{\pm 1}$ ($1 \leq j < n$)
satisfies

$$\hat \l(\w_j) = r^{\la\e_j,\l\ra}s^{\la\e_{j+1},\l\ra} \qquad \text{and} \qquad 
\hat \l(\w_j') = r^{\la\e_{j+1},\l\ra}s^{\la\e_j,\l\ra}. \tag 3.2$$ 
\m 
Let $M$ be a module for $\widetilde U = U_{r,s}(\fgl)$ of dimension $d < \infty$.
If $\K$ is algebraically closed ({\it which will be our
assumption throughout the remainder of this work}), then 

$$M = \bigoplus_{\chi} M_\chi$$

\n where each $\chi: \widetilde U^0 \rightarrow \K$ is an algebra homomorphism,
and  $M_\chi$ is the generalized eigenspace given by

$$M_\chi = \{m \in M \mid (a_i-\chi(a_i)\,1)^d m
= 0 =  (b_i-\chi(b_i)\,1)^d m, \ \ \text{ for all} \ i\}. \tag 3.3$$

\n  When $M_\chi \neq 0$ we say that $\chi$ is
a {\it weight} and $M_{\chi}$ is the corresponding {\it weight space}. 
(If $M$ decomposes into genuine eigenspaces relative
to $\widetilde U^0$ (resp. $U^0$), then we say that
$\widetilde U^0$ (resp. $U^0$) {\it acts semisimply on} $M$.)   
\m From relations (R2) and (R3) we deduce that 

$$\aligned
& e_j M_\chi \subseteq M_{\chi \cdot \widehat {\a_j}} \\
& f_j M_\chi \subseteq M_{\chi \cdot (\widehat {-\a_j})}\endaligned \tag 3.4$$

\n where $\widehat {\a_j}$ is as in (3.1), and $\chi \cdot \psi$ is the
homomorphism with values 
$(\chi \cdot \psi)(a_i) = \chi(a_i)\psi(a_i)$ \ 
and \ $(\chi \cdot \psi)(b_i) = \chi(b_i)\psi(b_i)$. 
In fact, if  $(a_i-\chi(a_i)\,1)^k m = 0$, then
$(a_i - \chi(a_i) r^{\la\e_i,\a_j\ra}\,1)^k e_j m = 0$,
and similarly for $b_i$ and for $f_j$.  This can be
used to show that the sum of the 
eigenspaces is a submodule of $M$, and so if $M$ is simple, this sum must
be $M$ itself.  Thus, in (3.3), we may replace the power $d$ by 1
whenever $M$ is simple, and $\widetilde U^0$ must act
semisimply in this case.    We also can see from (3.4) that for each simple $M$
there is a homomorphism $\chi$ so that all the weights of $M$ are
of the form $\chi \cdot \hat \zeta$, where $\zeta \in Q$. 

When all the weights of a module $M$ are of the form $\hat{\l}$, where
$\l\in\Lambda$, then for brevity we say that $M$ has weights in $\Lambda$.
Rather than writing $M_{\hat{\l}}$ for the weight space, we simplify the
notation by writing $M_{\l}$.  Note then (3.4) can be rewritten as
$e_jM_{\l}\subseteq M_{\l+\alpha_j}$ and $f_j M_{\l}\subseteq M_{\l-\alpha_j}$.
Any simple $\widetilde{U}$-module having one weight in $\Lambda$ has all its 
weights in $\Lambda$.

\m
 We would like to argue that
when $rs^{-1}$ is not a root of unity, the
elements $e_j$ and $f_j$ act nilpotently on any finite-dimensional
module.  For this we require the following result.  
\b
\proclaim {Proposition 3.5} Suppose $\hat \zeta = \hat \eta$,  where
$\zeta,\eta
\in Q$. If $rs^{-1}$ is not a root of unity, then $\zeta = \eta$.  \endproclaim
\m
 \demo{Proof} We will prove that when $\hat \zeta = \hat \eta$ as 
homomorphisms on the subalgebra $U^0$ generated by the $\w_i,\w_i'$, 
then $\zeta = \eta$, so the result holds for $U$ as well as for 
$\widetilde U$.  We may assume $\zeta =\sum_{i=1}^{n-1}\zeta_i \a_i$ 
and $\eta =\sum_{i=1}^{n-1}\eta_i \a_i$.  The condition that $\hat 
\zeta = \hat \eta$ gives the equations

$$\aligned
\hat \zeta(\w_i) & = r^{\la \e_i,\zeta\ra}s^{\la\e_{i+1},\zeta\ra} =
r^{\zeta_i-\zeta_{i-1}}s^{\zeta_{i+1}-\zeta_i} \\
& = \hat \eta(\w_i) =
r^{\eta_i-\eta_{i-1}}s^{\eta_{i+1}-\eta_i} \\
\hat \zeta(\w_i') & = r^{\la\e_{i+1},\zeta\ra}s^{\la\e_{i},\zeta\ra} =
r^{\zeta_{i+1}-\zeta_{i}}s^{\zeta_{i}-\zeta_{i-1}} \\
& = \hat \eta(\w_i') =
r^{\eta_{i+1}-\eta_{i}}s^{\eta_{i}-\eta_{i-1}}, \endaligned$$

\noindent where $\zeta_0 = \zeta_n = 0 = \eta_0 = \eta_n$.  Letting $\mu_i =
\zeta_i-\eta_i$, we may rewrite the above equations as
$$\align
& r^{\mu_i-\mu_{i-1}}s^{\mu_{i+1}-\mu_i}  = 1  \tag 3.6  \\
& r^{\mu_{i+1}-\mu_i}s^{\mu_i-\mu_{i-1}} = 1.  \tag 3.7 \endalign$$

\n Combining these we have
$$\align
& r^{\mu_{i+2}-\mu_{i+1}-\mu_i+\mu_{i-1}} = 1 \tag 3.8 \\
& s^{\mu_{i+2}-\mu_{i+1}-\mu_i+\mu_{i-1}} = 1 \tag 3.9 \endalign$$

\n for $i=1,\dots,n-2$.   Since we are assuming that $rs^{-1}$ is
not a root of unity, not both $r$ and $s$ are roots
of unity, so from these relations we see that

$$\mu_{i+2}-\mu_{i+1}-\mu_i+\mu_{i-1} = 0 \tag 3.10$$

\n for $i=1,\dots,n-2$.  We claim that the solution to the system of equations
given by (3.10) satisfies

$$\mu_{2k} =  k \mu_2  \qquad \text{and} \qquad \mu_{2k+1} = k\mu_2 + \mu_1. \tag 3.11$$ 

\n This is true for $\mu_3$ as $\mu_0 = 0$.  Moreover, $\mu_4 = \mu_3+\mu_2-\mu_1 = 2\mu_2$.
An easy induction proves the rest.  Now $\mu_n = 0$, and using that fact in (3.10) we
have
$$\mu_{n-1} = -\mu_{n-2} + \mu_{n-3}. \tag 3.12$$

If $n = 2m$ for some $m$,  then (3.11) and (3.12) give $\mu_2 = 0$.  From
(3.6), we have
$(rs^{-1})^{\mu_1} = 1$, and because $rs^{-1}$ is not a root of unity, this
says $\mu_1 = 0$.  The relations in (3.11) then show $\mu_i  = \zeta_i-\eta_i = 0$ for all
$i$.  Hence $\zeta = \eta$ when $n$ is even.  

Now if instead $n=2m+1$,  then (3.11) and (3.12) show that $\mu_1 = -m\mu_2$.  
The equations in (3.6) and (3.7) imply 
$$r^{-m\mu_2}s^{(m+1)\mu_2} = 1 \qquad \text{and} \qquad  r^{(m+1)\mu_2}s^{-m\mu_2} = 1, 
\tag 3.13$$

\n and hence that   $(rs)^{\mu_2} = 1.$    Then from (3.13) we see that
$s^{(2m+1)\mu_2} = 1 = r^{(2m+1)\mu_2}$.  As not both $r$ and $s$ are roots of unity, $\mu_2
= 0$. From this, the desired conclusion 
$\zeta = \eta$ follows.   \qed \enddemo
\b
\proclaim {Corollary 3.14} Let $M$ be a finite-dimensional module
for $U_{r,s}(\fsn)$ or for $U_{r,s}(\fgl)$.  If $rs^{-1}$ is not a root of unity,
then the elements $e_i, f_i$ ($1 \leq i < n)$ act nilpotently on $M$. \endproclaim
\m
\demo{Proof}  Because $M$ is a direct sum of its
weight spaces, it  suffices to argue that $e_i$ and $f_i$ act
nilpotently on each $M_{\chi}$.  As the weights $\widehat {k \a_i}$
for
$k=1,2,\dots$ are distinct by Proposition 3.5, and $e_i^k M_{\chi} \subseteq M_{\chi 
\cdot (\widehat {k\a_i})}$, it must be that some power of $e_i$ maps
$M_{\chi}$ to 0.  A similar argument applies to show that $f_i$ is
nilpotent also. 
\qed \enddemo
\m 
\b  
\head {\S 4. $R$-matrix and quantum Casimir operator} \endhead
\m 
Let $\Cal O$ denote the category of modules $M$ for $\widetilde U = U_{r,s}(\fgl)$
which satisfy the conditions:
\m \roster \item"{($\Cal O$1)}" $\widetilde U^0$ acts semisimply on 
$M$, and the set $\text{wt}(M)$ of weights of $M$ belongs to 
$\Lambda$: \quad $M = \bigoplus_{\l \in \text{wt}(M)} M_{\l}$, where 
$M_{\l}= \{m \in M \mid a_i.m = r^{\la \e_i,\l\ra}, \ \ b_i.m = s^{\la 
\e_i,\l \ra}$ for all $i\}$; \item"{($\Cal O$2)}" $\dim_{\K}€ M_{\l} < 
\infty$ for all $\l \in \text{wt}(M)$; \item"{($\Cal O$3)}" 
$\text{wt}(M) \subseteq \bigcup_{\mu \in F} (\mu - Q^+)$ for some 
finite set $F \subset \Lambda$.  \endroster \m \n The morphisms in $\Cal 
O$ are $\widetilde U$-module homomorphisms.  \m

Category $\Cal O$ is closed under tensor product. 
For any two modules $M$ and $M'$ in  
$\Cal O$, we construct a $\widetilde U$-module isomorphism $R_{M',M}: 
M' \ot M \rightarrow M \ot M'$, by the method used by Jantzen [Ja, 
Chap.  7] for the quantum groups $U_q(\g)$.  These isomorphisms work 
equally well for $U$-modules.  \m The map $R_{M',M}$ is the composite 
of three linear transformations $P$, $\widetilde f$, $\Theta$, which 
we now describe: \b \n (i) \ $P = P_{M',M}: M' \ot M \rightarrow M \ot 
M'$, \ \ $P(m' \ot m) = m\ot m'$.  \b \n (ii) \ $\widetilde f = 
\widetilde f_{M,M'}: M \ot M' \rightarrow M' \ot M$ is such that 
$\widetilde f(m \ot m') = f(\l,\mu)(m \ot m')$ when $m \in 
M_{\lambda}$ and $m' \in M'_{\mu}$, where the map $f: \Lambda \times 
\Lambda \rightarrow \K^\#$ is defined as follows.  \m Suppose that 
$\a_n = \e_n$ so that $\a_i + \a_{i+1} + \cdots + \a_n = \e_i$ for $i 
= 1,\dots,n$.  Let $\Lambda = \Z \a_1 \oplus \cdots \oplus \Z \a_n = 
\Z \e_1 \oplus \cdots \oplus \Z \e_n$, (the weight lattice of $\fgl$).  
If $\lambda=\sum_{i=1}^n\lambda_i\alpha_i$ is in the weight lattice 
$\Lambda$, we define $$\aligned \w_\lambda 
&=\w_1^{\lambda_1}\cdots\w_{n-1}^{\lambda_{n-1}}a_n^{\lambda_n}\\
\w_\lambda'&=(\w_1')^{\lambda_1}\cdots (\w_{n-1}')^{\lambda_{n-1}}
b_n^{\lambda_n},\\
\endaligned \tag 4.1$$
which agrees with (1.2) in case $\lambda\in Q$.  
If also $\mu=\sum_{i=1}^n\mu_i\alpha_i$ is in $\Lambda$, we define
$$ f(\lambda,\mu)=(\omega_\mu',\omega_\lambda)^{-1}. \tag 4.2 $$
The values of this bilinear form are given by (2.1) and (2.3).  It may
be checked that for all $\lambda, \mu, \nu\in \Lambda$ and $1\leq i,j <n$, 
the following hold:
$$\aligned
f(\lambda +\mu,\nu) & = f(\lambda,\nu)f(\mu,\nu)\\
f(\lambda,\mu+\nu) &= f(\lambda,\mu)f(\lambda,\nu)\\
f(\alpha_j,\mu)&=r^{-\la\epsilon_j,\mu\ra}s^{-\la\epsilon_{j+1},\mu\ra}\\
f(\lambda,\alpha_i) &= r^{\la\epsilon_{i+1},\lambda\ra}s^{\la\epsilon_i,\lambda\ra}.\\
\endaligned \tag 4.3 $$
 
We will need to compute 

$$f(\e_i,\e_j) = f(\a_i+\cdots + \a_n, \a_j+ \cdots + \a_n).$$

\n Supposing first that $1 \leq i,j < n$, by (4.3) we have

$$\aligned
f(\e_i,\e_j)&= (\w_{\alpha_j +\cdots +\alpha_n}',a_n)^{-1} \prod_{k=i}^{n-1}
f(\alpha_k,\alpha_j+\cdots +\alpha_n) \\
& = r\prod_{k = i}^{n-1} r^{-\la\e_k,\a_j+\cdots+\a_{n}\ra}
s^{-\la\e_{k+1},\a_j+\cdots+\a_{n}\ra} \\
& = r\prod_{k = i}^{n-1}
r^{-\la\e_k,\e_j\ra}s^{-\la\e_{k+1},\e_j\ra}  \\
& = \cases s^{-1} & \qquad \ \text{if} \ \ i < j \\
1 &  \qquad \ \text{if} \ \ i = j \\
r & \qquad \ \text{if} \ \ i > j.  \endcases
\endaligned$$

Now if $1 \leq i < n$, then

$$\aligned
f(\e_i,\e_n) & = f(\a_i+\cdots +\a_n,\a_n)=(b_n,\w_i)^{-1}\cdots
(b_n,\w_{n-1})^{-1}(b_n,a_n)^{-1}= s^{-1} \\
f(\e_n,\e_i) & = f(\a_n,\a_i+\cdots+\a_n)
= (\w_i,a_n)^{-1}\cdots (\w_{n-1},a_n)^{-1}(b_n,a_n)^{-1}= r \quad \quad 
\text{and} \\
f(\e_n,\e_n) & =(b_n,a_n)= 1. \endaligned $$
\m
\n As a result, the following holds:
\b
\proclaim{Lemma 4.4} For all $1 \leq i,j \leq n$,  we have
$$ 
f(\e_i,\e_j) = \cases s^{-1} & \qquad \ \text{if} \ \ i < j \\
1 &  \qquad \ \text{if} \ \ i = j \\
r & \qquad \ \text{if} \ \ i > j.  \endcases$$
\endproclaim 
\m
 \n (iii) \ Now we turn our attention to the construction of our 
final mapping.  Observe it is a consequence of (R2) and (R3) that the 
subalgebra $U^+$ of $\widetilde U$ (or of $U = U_{r,s}(\fsn)$) 
generated by $1$ and $e_i \ (1 \leq i < n)$ has the following 
decomposition

$$U^+ = \bigoplus_{\zeta \in Q^+} U^+_{\zeta}$$

\n where  
$$U^+_{\zeta} = \{z \in U^+ \mid a_i z = r^{\la\e_i,\zeta\ra}z a_i, \ \
b_i z = s^{\la\e_i,\zeta\ra} z b_{i}, \ \ (1 \leq i < n)\}.$$

\n The weight space $U^+_\zeta$ is spanned by all the monomials $e_{i_1} \ \cdots \
e_{i_\ell}$ such that $\a_{i_1} + \cdots + \a_{i_\ell} = \zeta$.
\m
Similarly, the subalgebra $U^-$ generated
by $1$ and the $f_i$'s has a decomposition $U^- = \bigoplus_{\zeta \in Q^+}
U^-_{-\zeta}$, and  the spaces $U^+_\zeta$ and $U^-_{-\zeta}$ are nondegenerately
paired. 
\m

 Since
$\Delta(e_i) = e_i \ot 1 + \w_i \ot e_i$, we have

$$\Delta(x) \in \sum_{0 \leq \nu \leq \zeta}  U^+_{\zeta-\nu}\omega_\nu \ot
U^+_\nu$$ 

\n for all $x\in U^+_{\zeta}$.
(In writing this, we are using the standard partial order on $Q$ in which
$\nu \leq \zeta$ if $\zeta-\nu \in Q^+$.)  For each $i$, there are elements
$p_i(x)$ and
$p_i'(x)
\in U^+_{\zeta-\a_i}$ such that

$$\aligned
\Delta(x) = x \ot 1 + \sum_{i=1}^{n-1} p_i(x)\w_i \ot e_i + \  \text{the rest} \\  
\Delta(x) = \w_\zeta \ot x  
+ \sum_{i=1}^{n-1}e_i \w_{\zeta-\a_i} \ot p_i'(x) + \  \text{the rest},
\endaligned \tag 4.5$$ 
where in each case ``the rest'' refers to terms involving products of more
than one $e_j$ in the second factor (respectively, in the first factor).
(Compare the expressions in Lemma 5.2 below.)
\b
\proclaim {Lemma 4.6} (Compare [Ja, Lemma 6.14, 6.17].)  For all $x \in U^+_\zeta$, \
$x'
\in U^+_{\zeta'}$, and $y \in U^-$, the following hold:
\roster
\item"{(i)}" $p_i(xx') = r^{\la\e_i,\zeta'\ra}s^{\la\e_{i+1},\zeta'\ra}p_i(x)x' + x p_i(x')$.
\s
\item"{(ii)}" $p_i'(xx') = p_i'(x)x' + r^{-\la\e_{i+1},\zeta\ra}s^{-\la\e_{i},\zeta\ra}x
p_i'(x)$.
\s
\item"{(iii)}" $(f_iy,x) = (f_i,e_i)(y,p_i'(x)) =  (s-r)^{-1}(y,p_i'(x))$.
\s
\item"{(iv)}" $(yf_i,x) = (f_i,e_i)(y,p_i(x)) =  (s-r)^{-1}(y,p_i(x))$. 
\s
\item"{(v)}" $f_i x - xf_i = (s-r)^{-1}\big(p_i(x)\w_i - \w_i' p_i'(x)\big)$. 
\endroster 
\endproclaim  
\m
 \demo{Proof} The proofs of (i) and (ii) amount to equating the 
expressions for $\Delta(xx') = \Delta(x)\Delta(x')$.  We demonstrate 
the second:

$$\aligned
\Delta(xx') & = \w_{\zeta+\zeta'} \ot xx'  
+ \sum_{i=1}^{n-1}e_i \w_{\zeta+\zeta'-\a_i} \ot p_i'(xx') + \  \text{the rest}\\
\Delta(x)\Delta(x') & = \Big(\w_\zeta \ot x  
+ \sum_{i=1}^{n-1}e_i \w_{\zeta-\a_i} \ot p_i'(x) + \  \text{the rest}\Big) \times \\
& \qquad \qquad  \Big(\w_{\zeta'} \ot x'  
+ \sum_{i=1}^{n-1}e_i \w_{\zeta'-\a_i} \ot p_i'(x') + \  \text{the rest}\Big)
\\ & =  \w_{\zeta+\zeta'} \ot xx' \\
& \quad \quad + \sum_{i=1}^{n-1} \Big( e_i
\omega_{\zeta+\zeta'-\a_i} \ot  p_i'(x)x'+ \w_\zeta e_i \w_{\zeta'-\a_i} \ot x
p_i'(x') \Big)  + \ \text{the rest} \\ & =  \w_{\zeta+\zeta'} \ot xx' \\
& \quad \quad  +
\sum_{i=1}^{n-1}  e_i
\omega_{\zeta+\zeta'-\a_i} \ot 
\left (p_i'(x)x'+ r^{-\la\e_{i+1},\zeta\ra}s^{-\la\e_{i},\zeta\ra}x p_i'(x')\right) +
\ \text{the rest}, \endaligned$$

\n (using Lemma 1.3).  Equating terms gives (ii).  

Now for (iii) we argue as follows using the second equation of (4.5):

$$\aligned
(f_iy, x) & = \sum (f_i, x_{(1)})(y, x_{(2)}) \\ 
& = (f_i,e_i \w_{\zeta-\a_i}) (y,p_i'(x)) \\
& = (1 \ot f_i + f_i \ot \w_i', \ \w_{\zeta-\a_i} \ot e_i) (y,p_i'(x)) \\
& = (f_i,e_i)(y,p_i'(x)). \endaligned$$

Let's begin the proof of (v) by observing that it is true if
$x = 1 \in U_0^+$ since $p_i(1) = p_i'(1) = 0$ for all $i$.  The
relation in (v) also holds for $x = e_j$, because $\Delta(e_j)
= e_j \ot 1 + w_j \ot e_j$ implies that $p_i(e_j) = \delta_{i,j} = p_i'(e_j)$.
We suppose the result is true for $x \in U^+_\zeta$ and
$x' \in U^+_{\zeta'}$ and prove it for $xx' \in U^+_{\zeta+\zeta'}$. Now

$$\aligned
f_ixx' - xx'f_i & = (f_ix-xf_i)x' - x(f_ix'-x'f_i) \\
& = \frac{1}{s-r}\Big(\big(p_i(x)\w_i - \w_i' p_i'(x)\big) x'
+x\big(p_i(x')\w_i - \w_i' p_i'(x')\big)\Big) \\ 
& = 
\frac{1}{s-r}\Bigg (\left(r^{\la\e_i,\zeta'\ra}s^{\la\e_{i+1},\zeta'\ra}p_i(x)x'
+ x\big(p_i(x')\right)\w_i \\
&\qquad \qquad \ \  - \w_i' \left(p_i'(x) x' +
r^{-\la\e_{i+1},\zeta\ra}s^{-\la\e_i,\zeta\ra} xp_i'(x')\right) \Bigg) \\
& = 
\frac{1}{s-r}\left(p_i(xx')\w_i - \w_i'p(xx') \right). \qed  \endaligned $$
\enddemo
Assuming for $y \in U^-_{-\zeta}$ that $p_i(y)$ and $p_i'(y)$ are defined by

$$\aligned 
\Delta(y) & = y \ot \w_\zeta' + \sum_{i=1}^n p_i(y) \ot f_i \w'_{\zeta -\a_i}+
\ \text{the rest} \\
\Delta(y) & = 1 \ot y + \sum_{i=1}^n f_i \ot p_i'(y) \w'_i   + \
\text{the rest},\endaligned \tag 4.7$$

\n the same type of argument produces this analogue of Lemma 4.6:
 
\b
\proclaim {Lemma 4.8} For all $y \in U^-_{-\zeta}$, \
$y'
\in U^-_{-\zeta'}$, and $x \in U^+$, the following hold:
\roster
\item"{(i)}" $p_i(yy') = p_i(y)y' + r^{\la\e_i,\zeta'\ra}s^{\la\e_{i+1},\zeta'\ra}y p_i(y')$.
\s
\item"{(ii)}" $p_i'(yy') = r^{-\la\e_{i+1},\zeta'\ra}s^{-\la\e_{i},\zeta'\ra}p_i'(y)y' + y
p_i'(y)$.
\s
\item"{(iii)}" $(y,e_ix) = (f_i,e_i)(p_i(y),x)=
 (s-r)^{-1}(p_i(y),(x))$.
\s
\item"{(iv)}" $(y,xe_i) = (f_i,e_i)(p_i'(y),x) =
 (s-r)^{-1}(p_i'(y),x)$. 
\s
\item"{(v)}" $e_i y - ye_i = (r-s)^{-1}\big(w_ip_i(y)- p_i'(y)\w_i' \big)$. 
\endroster 
\endproclaim 
\b  

Because the spaces $U^+_\zeta$ and $U^-_{-\zeta}$ are nondegenerately
paired,  we may select a basis $\{u_k^\zeta\}_{k=1}^{d_\zeta}$,  
$(d_\zeta = \dim_\K U^+_\zeta)$, for
$U^+_\zeta$ and a dual
basis
$\{v_k^\zeta\}_{k=1}^{d_\zeta}$ for $U^-_{-\zeta}$. 
Then for each $x \in U^+_\zeta$ and $y\in U^-_{-\zeta}$ we have

$$x = \sum_{k=1}^{d_\zeta} (v_k^\zeta,x) u_k^\zeta \ \ \ \text{and} \ \ \
y=\sum_{k=1}^{d_{\zeta}}(y,u_k^{\zeta})v_k^{\zeta}. \tag 4.9$$
\m
For $\zeta \in 
Q^+ = \bigoplus_{i=1}^{n-1} \Z_{\geq 0} \a_i$, we define  

$$\Theta_{\zeta} = \sum_{k = 1}^{d_\zeta} v_k^\zeta \ot u_k^{\zeta}.$$

\n  Set $\Theta_\zeta = 0$ if $\zeta \not \in Q^+$. 
\b
\proclaim {Lemma 4.10} $\ii$  $(a_i \ot a_i) \Theta_\zeta =
\Theta_\zeta (a_i \ot a_i)$;
\roster
\item"{(ii)}" $(b_i \ot b_i) \Theta_\zeta =
\Theta_\zeta (b_i \ot b_i)$; \s

\item"{(iii)}" $(e_i \ot 1) \Theta_{\zeta} + (\w_i \ot e_i) \Theta_{\zeta-\a_i}
=  \Theta_{\zeta}(e_i \ot 1) + \Theta_{\zeta-\a_i}(\w_i' \ot e_i)$; \s
\item"{(iv)}" $(1 \ot f_i)\Theta_\zeta + (f_i \ot \w_i')\Theta_{\zeta-\a_i}
= \Theta_{\zeta}(1 \ot f_i) + \Theta_{\zeta-\a_i}(f_i \ot \w_i);$
\endroster
\n for $1 \leq i < n$.  
\endproclaim
\m
\demo{Proof} The first two are easy to check.  We demonstrate (iv)
and leave (iii) as an exercise. The calculation below will use (iii)-(v)
of Lemma 4.6 and  (4.9). 

$$\allowdisplaybreaks 
\aligned  
(1 \ot f_i) \Theta_\zeta  - \Theta_\zeta & (1 \ot f_i) \\
& = \sum_{k} v_k^\zeta \ot (f_i u_k^\zeta - u_k^\zeta f_i) \\
& = \frac{1}{s-r} \sum_k v_k^\zeta \ot \left (p_i(u_k^\zeta)\w_i
-\w_i'p_i'(u_k^\zeta)\right) \\
& = \frac{1}{s-r} \sum_k v_k^\zeta \ot \left (\sum_j \big (v_j^{\zeta-\a_i},\,
p_i(u_k^\zeta)\big) u_j^{\zeta-\a_j}\right)\w_i \\
& \qquad \qquad 
- \frac{1}{s-r} \sum_k v_k^\zeta \ot \w_i'\left (\sum_j \big (v_j^{\zeta-\a_i},\,
p_i'(u_k^\zeta)\big) u_j^{\zeta-\a_j}\right) \\
& = \sum_k v_k^\zeta \ot \left( \sum_j
(v_j^{\zeta-\a_i}f_i,\,u_k^\zeta) u_j^{\zeta-\a_i}\right)\w_i \\
& \qquad \qquad  - \sum_k v_k^\zeta \ot
\w_i'\left (\sum_j \big (f_iv_j^{\zeta-\a_i},\,
 u_k^\zeta\big) u_j^{\zeta-\a_j}\right)\\ 
 & = \sum_j \left ( \sum_k (v_j^{\zeta-\a_i}f_i,\,u_k^\zeta) 
 v_k^\zeta \right) \ot u_j^{\zeta-\a_i}\w_i \\
& \qquad \qquad 
\ - \ \sum_j \left(\sum_k (f_iv_j^{\zeta-\a_i},\,
 u_k^\zeta\big) v_k^\zeta \right)  \ot \w_i' u_j^{\zeta-\a_j} \endaligned$$
 
$$\aligned
 \hskip 1.5 truein 
& = \sum_j v_j^{\zeta-\a_i}f_i \ot u_j^{\zeta-\a_i}\w_i 
\ - \ \sum_j f_iv_j^{\zeta-\a_i} \ot \w_i' u_j^{\zeta-\a_j} \\
& = \Theta_{\zeta-\a_i}(f_i \ot \w_i) - (f_i\ot \w_i')\Theta_{\zeta-\a_i}. 
\endaligned$$

\n Consequently, (iv) holds. \qed 
\enddemo 
\b
We now define 

$$\Theta = \sum_{\zeta \in Q^+} \Theta_\zeta$$

\n  (which we can think of as living in the completion of
$\widetilde{U} \ot \widetilde{U}$  where infinite sums are allowed). 
For fixed $\widetilde U$-modules $M$ and $M'$ in $\Cal O$, we may apply 
$\Theta$ to their tensor product:

$$\Theta = \Theta_{M,M'}: M \ot M' \rightarrow M \ot M'. $$

\n  Note that $\Theta_\zeta: M_{\l} \ot M'_{\mu} \rightarrow 
M_{\l - \zeta} \ot M'_{\mu+\zeta}$ for all $\l,\mu \in \Lambda$, 
and
because of condition
$({\Cal O}3)$, there are only finitely many $\zeta \in Q^+$ such that 
$M'_{\mu+\zeta} \neq 0$.  Hence, this is a well-defined 
linear transformation on $M \ot M'$.  
 \m We can choose countable bases of weight vectors for both $M$ and $M'$ and 
 their tensor products as a basis for $M \ot M'$.  Then ordering this 
 basis appropriately shows that each $\Theta_\zeta$ with $\zeta > 0$ 
 has a strictly upper triangular matrix.  Because $\Theta_0 = 1 \ot 1$ 
 acts as the identity transformation on $M \ot M'$, \ $\Theta_{M,M'}$ 
 is an invertible transformation.

\b \proclaim 
{Theorem 4.11} Let $M$ and $M'$ be modules in $\Cal O$.  Then the map

$$\Theta \circ \widetilde f \circ P: M' \ot M \rightarrow M \ot M'$$

\n is an isomorphism of $\widetilde U$-modules.  \endproclaim 
\m
 \demo{Proof} Since each of the maps is invertible, once we show 
that $\Theta \circ \widetilde f \circ P$ is a $\widetilde U$-module 
homomorphism, we will be done.  The proof amounts to verifying that

$$\Delta(a)(\Theta \circ \widetilde f \circ P)(m' \ot m) = (\Theta \circ \widetilde f
\circ P)
\Delta(a)(m' \ot m) \tag 4.12$$

\n holds for all $a \in \widetilde U$, $m \in M_{\l}$ and
$m' \in M_{\mu}'$.   Because $\Delta$ is an algebra homomorphism, it suffices to
check (4.12) on the generators $e_i,f_i, a_i, b_i$. We will present the computation
just for $a = e_i$.  In this case,  the right side of (4.12) becomes

$$\aligned
(\Theta \circ \widetilde f
\circ P)
\Delta(e_i)(m' \ot m) & = (\Theta \circ \widetilde f
\circ P)(e_i m' \ot m + \w_i m' \ot e_i m) \\
& = (\Theta \circ \widetilde f)(m \ot e_i m' + e_i m \ot \w_i m') \\
& = f(\l,\mu+\a_i) \Theta(m \ot e_i m') + 
f(\l + \a_i,\mu) \Theta(e_i m \ot \w_i m') \\ 
& = f(\l,\mu+\a_i)\Big(\sum_{\zeta} \Theta_\zeta \Big) (1\ot e_i) (m  \ot m') \\
& \qquad \qquad 
+ f(\l + \a_i,\mu)\Big(\sum_{\zeta}\Theta_\zeta \Big)(e_i \ot \w_i)(m \ot m') \\
\endaligned $$ 

Now let's compute the left side using (iii) of Lemma 4.10:

$$\aligned
\Delta(e_i)(\Theta \circ \widetilde f
\circ P) (m' \ot m) & = f(\l,\mu)\Delta(e_i)\Theta(m \ot m')
\\ & = f(\l,\mu) (e_i \ot 1) \Big(\sum_{\zeta} \Theta_\zeta \Big)(m \ot m') \\
& \qquad \qquad  + f(\l,\mu)(\w_i \ot e_i) \Big(\sum_{\zeta} \Theta_{\zeta-\a_i}
\Big)(m
\ot m')
\\ & = f(\l,\mu) \Big(\sum_\zeta \Theta_\zeta\Big)(e_im \ot m')   \\
& \qquad \qquad + f(\l,\mu) \Big(\sum_{\zeta} \Theta_{\zeta-\a_i}\Big)
(\w_i'm \ot e_i m') \\
& = f(\l,\mu)r^{-\la\e_i,\mu\ra}s^{-\la\e_{i+1},\mu\ra} \Big(\sum_\zeta
\Theta_\zeta\Big)(e_im
\ot
\w_im')  \\
& \qquad \qquad + f(\l,\mu)r^{\la\e_{i+1},\l\ra}s^{\la\e_i,\l\ra} \Big(\sum_{\zeta}
\Theta_{\zeta-\a_i}\Big) (m \ot e_i m'). 
\endaligned $$ 

\n This expression can be seen to equal the previous one by (4.3).
(Note in this computation we have
made liberal use of the fact that $\sum_\zeta \Theta_\zeta = 
\sum_\zeta \Theta_{\zeta -\a_i}$ because of our convention that
$\Theta_\eta = 0$ whenever $\eta \not \in 	Q^+$.)   \qed  \enddemo
\m
\subhead Quantum Casimir operator \endsubhead \m

In this subsection, we construct a quantum Casimir operator which commutes 
with the action of $\widetilde U$ on any $\widetilde
U$-module in $\Cal O$.  \m

 Again letting 
$\{u_k^{\zeta}\}_{k=1}^{d_{\zeta}}$ and $\{v_k^{\zeta}\}_{k=1}^ 
{d_{\zeta}}$ be dual bases for $U_{\zeta}^+$ and $U^-_{-\zeta}$ 
respectively ($d_{\zeta}=\dim_{\K}€U_{\zeta}^+$), define 
$$\Omega=\sum_{\zeta\in 
Q^+}\sum_{k=1}^{d_{\zeta}}S(v^{\zeta}_k)u^{\zeta}_k,
\tag 4.15$$
where $S$ denotes the antipode. Note that by (3.4), $\Omega$ preserves the weight spaces
of any $M \in \Cal O$.  \b \proclaim{Lemma 4.16} For $M \in 
\Cal O$, assume  $m\in M_{\l}$.  Then \roster 
\item"{(i)}" $\Omega e_i .  m = (rs^{-1})^{-\langle 
\alpha_i,\l+\alpha_i\rangle } e_i \Omega .  m$,

\item"{(ii)}" $\Omega f_i .  m = (rs^{-1})^{\langle 
\alpha_i,\l\rangle}f_i\Omega .  m$.  \endroster \endproclaim 
\m 
\demo{Proof} Apply $\frak{m} \circ (S\otimes 1)$ to Lemma 4.10 (iii), 
where $\frak{m}$ is the multiplication map.  As $S$ is an algebra 
anti-automorphism, this yields $$\aligned 
-\sum_{k=1}^{d_{\zeta}}S(v^{\zeta}_k)\omega_i^{-1}e_iu^{\zeta}_k & + 
\sum_{k=1}^{d_{\zeta 
-\alpha_i}}S(v^{\zeta-\alpha_i}_k)\omega_i^{-1}e_iu^{\zeta 
-\alpha_i}_k \\
& \quad 
=-\sum_{k=1}^{d_{\zeta}}\omega_i^{-1}e_iS(v^{\zeta}_k)u^{\zeta}_k + 
\sum_{k=1}^{d_{\zeta-\alpha_i}}(\omega_i')^{-1}S(v_k^{\zeta 
-\alpha_i})
u^{\zeta -\alpha_i}_k e_i.  \endaligned$$
Now act on $m \in M_{\lambda}€$ with the result, and sum over all 
$\zeta\in Q^+$.  The two sums on the left side cancel, while the right 
side produces
$$\omega_i^{-1}e_i\Omega . m = (\omega_i')^{-1}\Omega e_i . m.$$
By definition, $\Omega$ preserves weight spaces, and so $e_i\Omega . m, \
\Omega e_i . m\in M_{\l + \alpha_i}$.  Therefore we have
$$r^{-\langle\epsilon_i,\l +\alpha_i\rangle}s^{-\langle\epsilon_{i+1},\l
+\alpha_i\rangle} e_i\Omega . m =r^{-\langle\epsilon_{i+1},\l +\alpha_i\rangle}
s^{-\langle\epsilon_{i},\l+\alpha_i\rangle}\Omega e_i . m,$$
which is equivalent to (i).

The proof of (ii) is virtually identical, and uses Lemma 4.10 (iv).  \qed 
\enddemo
\b

Now we introduce a certain function $g:\Lambda\rightarrow 
\K^\#$. If $\rho$ denotes the half sum of the positive roots,
then $2\rho  
 = \sum_{j=1}^n (n+1-2j)\e_j \in \Lambda$.  For 
$\lambda \in \Lambda$, set

$$g(\lambda) = (rs^{-1})^{\half \la \lambda + 2 \rho, \lambda 
\ra}. \tag 4.17$$

\n Because $\la \rho, \a_i \ra = 1$ for all $i=1,\dots,n-1$, 
it is straightforward to verify that

$$ g(\l +\alpha_i) = (rs^{-1})^{\langle\alpha_i,\l+\alpha_i\rangle} g(\l) 
\tag 4.18 $$
for all $\l\in\Lambda, \ i \in\{1,\ldots,n-1\}$.  
\m
For $M \in \Cal O$, define the linear operator $\Xi :M\rightarrow M$ by
$$ \Xi (m) = g(\l)m \tag 4.19$$
for all $m\in M_{\l}, \ \l\in\Lambda$.
\b  \proclaim{Theorem 4.20} The operator $\Omega\Xi:M\rightarrow M$
commutes with the action of $\widetilde U$ on 
any module $M \in \Cal O$.  \endproclaim

\m 
\demo{Proof} As $\Omega\Xi$ preserves the weight spaces of $M$, it 
commutes with the action of $a_i,b_i \ (1\leq i\leq n)$.  It remains 
to show that $\Omega\Xi$ commutes with $e_i,f_i \ (1\leq i<n)$.  Let 
$m\in M_{\l}$.  By Lemma 4.16 (i) and (4.18), we have 

$$\aligned 
\Omega\Xi(e_i .  m) & = g(\l +\alpha_i)\Omega e_i .  m\\
 &= (rs^{-1})^{\langle\alpha_i,\l +\alpha_i\rangle}g(\l)\Omega e_i .  
 m\\
&=g(\l)e_i\Omega . m = e_i\Omega \Xi (m).
\endaligned $$
The calculation for $f_i$ is similar. \qed\enddemo

 \b \m \head {\S 5.  Yang-Baxter equation 
and hexagon identities} \endhead \m  For pairs $M$, $M'$ of 
$\widetilde U$-modules in
category $\Cal O$, we will show first that the maps  
${R}_{M,M'}=\Theta\circ \widetilde{f}\circ P:M\ot M'\rightarrow M'\ot 
M$ satisfy the quantum Yang-Baxter equation.  That is, given three 
$\widetilde U$-modules $M$, $M'$, $M''$ in $\Cal O$, we have $R_{12}\circ 
R_{23}\circ R_{12} = R_{23}\circ R_{12}\circ R_{23}$ as maps from 
$M\ot M'\ot M''$ to $M''\ot M'\ot M$.  This abbreviated notation is 
standard, for example $R_{12}$ is an application of $R_{M,M'}$ to the 
first two of three factors and the identity map on the third factor.  
\m We will need the following inner product relation.  If $x\in 
U^+_{\gamma}$, $y\in U^-_{-\gamma}$, and $\zeta, \eta \in Q$, then

$$ (y\omega_{\zeta}', x\omega_{\eta}) = (y,x)(\omega_{\zeta}',\omega_{\eta}).
\tag 5.1$$

\n To derive this, we apply (4.5) and (4.7), keeping in mind that we need 
to take the {\it opposite} coproduct in the first position (equivalently,
reverse the order of the factors in the second position):
$$\aligned
(y\omega_{\zeta}',x\omega_{\eta}) &= (y\otimes\omega_{\zeta}',\Delta(x)
\Delta(\omega_{\eta}))\\
&= (y,x\omega_{\eta})(\omega_{\zeta}',\omega_{\eta})\\
&= (\Delta(y),\omega_{\eta}\ot x)(\omega_{\zeta}',\omega_{\eta})\\
&= (y,x)(\omega_{\zeta}',\omega_{\eta}). \endaligned $$
\b
\proclaim{Lemma 5.2} Let $x\in U^+_{\gamma}$ and $y\in U^-_{-\gamma}$. Then
\roster
\item"{(i)}"  $\Delta(x) = \sum_{0\leq \zeta\leq \gamma}\sum_{i,j}(v_i^{\gamma-\zeta}
v_j^{\zeta},x)u_i^{\gamma -\zeta}\omega_{\zeta}\otimes u_j^{\zeta}$,
\smallskip
\item"{(ii)}"  $\Delta(y) = \sum_{0\leq\zeta\leq\gamma}\sum_{i,j} (y,u_i^{\gamma-\zeta}
u_j^{\zeta})v_j^{\zeta}\otimes v_i^{\gamma-\zeta}\omega_{\zeta}'.$
\endroster
\endproclaim
\m
\demo{Proof} As $x\in U^+_{\gamma}$, we have $\Delta(x) \in \sum_{0\leq
\zeta\leq\gamma}U^+_{\gamma-\zeta}\omega_{\zeta}\otimes U^+_{\zeta}$.
Let $c_{ij}^{\zeta}\in \K$ be such that
$$ \Delta(x)=\sum_{\zeta,i,j}c_{ij}^{\zeta} u_i^{\gamma-\zeta}\omega_{\zeta}
\otimes u_j^{\zeta}.$$
Then for all $k,\ell,$ and  $\nu$, we see from (5.1) that
$$\aligned
(v_k^{\gamma-\nu}v_{\ell}^{\nu},x) &= (v_k^{\gamma-\nu}\otimes v_{\ell}^{\nu},
\Delta(x))\\
&= \sum_{\zeta,i,j}c_{i,j}^{\zeta}(v_k^{\gamma-\nu},u_i^{\gamma-\zeta}\omega
_{\zeta})(v_{\ell}^{\nu},u_j^{\zeta})= c_{k\ell}^{\nu},\endaligned $$

\n which proves (i).  The argument for (ii) is similar. \qed \enddemo
\b
Letting $\Theta^{\op}=\sum_{\gamma\in Q^+}\sum_i u_i^{\gamma}\ot v_i^
{\gamma}$, \
$\Theta_{12} =\sum_{\gamma\in Q^+}\sum_i v_i^{\gamma}\ot u_i^{\gamma}\ot 1$, \ 
$\Theta_{ij}^f = \Theta_{ij}\circ \widetilde{f}_{ij}$, and defining
the other expressions in a like manner, we have the following identities
for operators on $M\ot M'\ot M''$.
\b
\proclaim{Lemma 5.3} $\ii$\, $(\Delta\ot 1)(\Theta^{\op}) \circ \widetilde{f}_{31}
\circ \widetilde{f}_{32} = \Theta^f_{31}\circ \Theta^f_{32}.$
\roster
\item"{(ii)}" $\widetilde{f}_{31}\circ\widetilde{f}_{32}\circ\Theta_{12} =\Theta_{12}\circ
\widetilde{f}_{31}\circ\widetilde{f}_{32}.$
\endroster
\endproclaim
\m
 \demo{Proof} Let $m\in M_{\lambda}$, $m'\in M'_{\mu}$, and $m''\in 
M''_{\nu}$.  Then by Lemma 5.2(i) and (4.9), the left side of (i) 
applied to $m\ot m'\ot m''$ is $$\aligned (\Delta\ot 
1)(\Theta^{\op})\circ \widetilde{f}_{31}\circ\widetilde{f}_{32} &= 
f(\nu,\mu)f(\nu,\lambda)(\Delta\ot 1)\big (\sum_{\gamma, 
k}u_k^{\gamma}\ot v_k^{\gamma}\big )\\
&= f(\nu,\mu)f(\nu,\lambda)\sum_{\gamma,k}\sum_{\zeta,i,j}(v_i^{\gamma-\zeta}
v_j^{\zeta},u_k^{\gamma})u_i^{\gamma-\zeta}\omega_{\zeta}\ot u_j^{\zeta}
\ot v_k^{\gamma}\\
&= f(\nu,\mu)f(\nu,\lambda)\sum_{\gamma,\zeta,i,j}u_i^{\gamma-\zeta}\omega
_{\zeta}\ot u_j^{\zeta}\ot(\sum_k \big (v_i^{\gamma-\zeta}v_j^{\zeta},u_k^{\gamma})
v_k^{\gamma}\big)\\
&=  f(\nu,\mu)f(\nu,\lambda)\sum_{\gamma,\zeta,i,j}u_i^{\gamma-\zeta}\omega
_{\zeta}\ot u_j^{\zeta}\ot v_i^{\gamma-\zeta}v_j^{\zeta}.\endaligned $$
On the other hand,
$$\aligned
\Theta_{31}^f\circ \Theta_{32}^f(m\ot m'\ot m'') &= f(\nu,\mu)
\sum_{\eta,\zeta,i,j}f(\nu-\zeta,\lambda) u_i^{\eta}m\ot 
u_j^{\zeta}m'\ot v_i^{\eta}v_j^{\zeta}m''\\
&= f(\nu,\mu)f(\nu,\lambda)\sum_{\eta,\zeta,i,j}f(-\zeta,\lambda)u_i^{\eta}
\ot u_j^{\zeta}\ot v_i^{\eta}v_j^{\zeta}(m\ot m'\ot m'').\endaligned $$
Changing variables in the first expression above to $\eta=\gamma-\zeta$,
and noticing that $\omega_{\zeta}.m=f(-\zeta,\lambda)m$, we obtain the
second expression, proving (i).

Identity (ii) results from a simple calculation using (4.3).
 \qed  \enddemo
\b
We are now ready to verify the quantum Yang-Baxter equation.
\b
\proclaim{Theorem 5.4} (Compare [Ja, \S 7.6].)  
$R_{12}\circ R_{23}\circ R_{12}=R_{23}\circ
R_{12}\circ R_{23}$ as maps from $M\ot M'\ot M''$ to $M''\ot M'\ot M$.
\endproclaim
\m
 \demo{Proof} Note that $P_{\sigma}\circ 
\Theta_{ij}^f=\Theta^f_{\sigma(i)\sigma(j)}\circ P_{\sigma}$ for all 
permutations $\sigma$, and that the $\widetilde{f}_{ij}$ commute with 
one another.  Applying Lemma 5.3 and Theorem 4.11, we have $$\aligned 
R_{12}\circ R_{23}\circ R_{12} &= P_{12}\circ 
P_{23}\circ\Theta^f_{31}\circ \Theta^f_{32}\circ R_{12}\\
&= P_{12}\circ P_{23} \circ(\Delta\ot 1)(\Theta^{\op})\circ \widetilde{f}_{31}
\circ \widetilde{f}_{32}\circ \Theta_{12}\circ\widetilde{f}_{12}\circ P_{12}\\
&= P_{12}\circ P_{23}\circ (\Delta\ot 1)(\Theta^{\op})\circ \Theta_{12}\circ
\widetilde{f}_{31}\circ\widetilde{f}_{32}\circ\widetilde{f}_{12}\circ P_{12}\\
&= P_{12}\circ P_{23}\circ (\Delta\ot 1)(\Theta^{\op})\circ
\Theta_{12}\circ\widetilde{f}
_{12}\circ P_{12}\circ \widetilde{f}_{32}\circ \widetilde{f}_{31}\\
&= P_{12}\circ P_{23}\circ (\Delta\ot 1)(\Theta^{\op})\circ R_{12}\circ
\widetilde{f}_{32}\circ \widetilde{f}_{31}\\
&= P_{12}\circ P_{23}\circ R_{12}\circ (\Delta\ot 1)(\Theta^{\op})\circ
\widetilde{f}_{32}\circ \widetilde{f}_{31}\\
&= P_{12}\circ P_{23}\circ \Theta^f_{12}\circ P_{12}\circ \Theta^f_{31}\circ
\Theta^f_{32}\\
&= \Theta_{23}^f\circ P_{12}\circ P_{23}\circ P_{12}\circ \Theta^f_{31}
\circ \Theta^f_{32}\\
&= \Theta_{23}^f\circ P_{23}\circ P_{12}\circ P_{23}\circ \Theta^f_{31}
\circ \Theta^f_{32}\\
&= \Theta_{23}^f\circ P_{23}\circ \Theta^f_{12}\circ P_{12}\circ \Theta_{23}^f
\circ P_{23}\\
&= R_{23}\circ R_{12}\circ R_{23}. \qed \endaligned $$ \enddemo
\b
Next we will verify the hexagon identities.  For this we require two additional 
lemmas regarding operators on $M\ot M'\ot M''$.
\b
\proclaim{Lemma 5.5} $(\Delta\ot 1)(\Theta_{\gamma}) = \sum_{0\leq \zeta
\leq \gamma}(\Theta_{\gamma-\zeta})_{23}(\Theta_{\zeta})_{13}(1\ot \omega'
_{\zeta}\ot 1),$ and

\n $(1\ot \Delta)(\Theta_{\gamma}) = \sum_{0\leq \zeta\leq \gamma}(\Theta
_{\gamma-\zeta})_{12}(\Theta_{\zeta})_{13}(1\ot\omega_{\zeta}\ot 1).$
\endproclaim
\m

\demo{Proof} By the definition of $\Theta_{\gamma}$, Lemma 5.2(ii), and (4.9),
we have
$$\aligned
(\Delta\ot 1)(\Theta_{\gamma}) & =\sum_k\Delta(v_k^{\gamma})\ot u_k^{\gamma}\\
 &= \sum_k\sum_{\zeta,i,j}(v_k^{\gamma},u_i^{\gamma-\zeta}u_j^{\zeta})v_j
^{\zeta}\ot v_i^{\gamma-\zeta}\omega_{\zeta}'\ot u_k^{\gamma}\\
&= \sum_{\zeta,i,j}v_j^{\zeta}\ot v_i^{\gamma-\zeta}\omega_{\zeta}'\ot
(\sum_k (v_k^{\gamma},u_i^{\gamma-\zeta}u_j^{\zeta})u_k^{\gamma})\\
&= \sum_{\zeta,i,j}v_j^{\zeta}\otimes v_i^{\gamma-\zeta}\omega_{\zeta}'\ot
u_i^{\gamma-\zeta}u_j^{\zeta}\\
&= \sum_{0\leq \zeta\leq \gamma}(\Theta_{\gamma-\zeta})_{23}(\Theta_{\zeta})
_{13}(1\ot \omega_{\zeta}'\ot 1).\endaligned $$

The second identity may be checked in just the same way.\qed \enddemo
\b
\proclaim{Lemma 5.6}  $\widetilde{f}_{12}\circ (\Theta_{\eta})_{13} =(\Theta
_{\eta})_{13}\circ (1\ot \omega_{\eta}\ot 1)\circ \widetilde{f}_{12},$ and

\n $\widetilde{f}_{23}\circ (\Theta_{\eta})_{13} = (\Theta_{\eta})_{13}\circ
(1\ot \omega_{\eta}'\ot 1)\circ \widetilde{f}_{23}.$\endproclaim
\m
 \demo{Proof} Let $m\in M_{\lambda}$, $m'\in M'_{\mu}$, $m''\in 
M''_{\nu}$.  Then $$\aligned \widetilde{f}_{12}\circ 
(\Theta_{\eta})_{13}(m\ot m'\ot m'') &= f(\lambda-\eta,\mu)\sum_i 
v_i^{\eta}m\ot m'\ot u_i^{\eta}m''\\
&= f(\lambda,\mu)f(-\eta,\mu)\sum_iv_i^{\eta}m\ot m'\ot u^{\eta}_i m''\\
&= f(\lambda,\mu)\sum_iv_i^{\eta}m\ot \omega_{\eta}m'\ot u_i^{\eta}m''\\
&= (\Theta_{\eta})_{13}\circ (1\ot \omega_{\eta}\ot 1)\circ\widetilde{f}_{12}(m
\ot m'\ot m'').\endaligned $$

The second identity can be shown using $\omega_{\eta}'m'=f(\mu,\eta)m'$.
\qed\enddemo
\b
We continue with our assumption that $M, M'$, and $M''$ 
are $\widetilde{U}$-modules in $\Cal O$.  To verify the hexagon 
identities, let $\widetilde{f}'$ denote the transformation on $M''\ot 
M\ot M'$ taking $m''\ot m\ot m'\in M''_{\nu}\ot M_{\lambda}\ot 
M_{\mu}'$ to $f(\nu,\lambda +\mu)m''\ot m\ot m'$.  Let 
$\widetilde{f}''$ be the transformation on $M'\ot M''\ot M$ taking 
$m'\ot m''\ot m \in M'_{\mu}\ot M''_{\nu}\ot M_{\lambda}$ to 
$f(\mu+\nu,\lambda)m'\ot m''\ot m$.  As in [Ja, Thm.  3.18], the 
hexagon identities are equivalent to $R_{12}\circ R_{23} = (1\ot 
\Delta)(\Theta)\circ \widetilde{f}'\circ P_{12} \circ P_{23}$ as maps 
from $M\ot (M'\ot M'')\rightarrow (M''\ot M)\ot M'$, and $R_{23}\circ 
R_{12}=(\Delta\ot 1)(\Theta)\circ\widetilde{f}''\circ P_{23} \circ 
P_{12}$ as maps from $(M\ot M')\ot M''$ to $M'\ot(M''\ot M)$.  \b 
\proclaim{Theorem 5.7} The hexagon identities hold, that is, \roster 
\item"{(i)}" $R_{12}\circ R_{23} = (1\ot \Delta)(\Theta)\circ 
\widetilde{f}'\circ P_{12} \circ P_{23}$, and \smallskip \item"{(ii)}" 
$R_{23}\circ R_{12}=(\Delta\ot 1)(\Theta)\circ\widetilde{f}''\circ 
P_{23} \circ P_{12}$.  \endroster \endproclaim \m

\demo{Proof}  Let $m\ot m'\ot m''\in M_{\lambda}\ot M'_{\mu}\ot M''_{\nu}$.
By Lemma 5.5, the right side of (i) applied to $m\ot m'\ot m''$ gives
$$\aligned
(1\ot \Delta)(\Theta)\circ & \widetilde{f}'\circ P_{12}\circ P_{23}(m\ot m'\ot m'') \\
& = f(\nu,\lambda +\mu) \sum_{\gamma \in Q^+}\sum _{0\leq \zeta\leq \gamma}
(\Theta_{\gamma-\zeta})_{12}(\Theta_{\zeta})_{13}(1\ot \omega_{\zeta}\ot 1)
(m''\ot m\ot m'). \endaligned$$
On the other hand,  by Lemma 5.6, the left side of (i) can be seen to equal 
$$\aligned
\Theta_{12}\circ\widetilde{f}_{12}\circ P_{12}\circ & \Theta_{23}\circ  \widetilde{f}_{23}
\circ P_{23}(m\ot m'\ot m'') = \Theta_{12}\circ \widetilde{f}_{12}\circ
\Theta_{13}\circ\widetilde{f}_{13}(m''\ot m\ot m')\\
&= \sum_{\eta \in Q^+} \Theta_{12}\circ (\Theta_{\eta})_{13}\circ (1\ot\omega
_{\eta}\ot 1)\circ\widetilde{f}_{12}\circ\widetilde{f}_{13}(m''\ot m\ot m')\\
&= f(\nu,\lambda)f(\nu,\mu)\sum_{\eta \in Q^+}\Theta_{12}\circ(\Theta_{\eta})
_{13}\circ (1\ot \omega_{\eta}\ot 1)(m''\ot m\ot m'). \endaligned $$

\n Then because $f(\nu,\lambda +\mu) = f(\nu,\lambda)f(\nu,\mu)$, 
a change of variables
shows that this is equal to the right side of (i).
The proof of (ii) is similar.\qed\enddemo
\b
\n {\bf Remark 5.8.}  $\Cal O$ is a braided monoidal category with 
braiding $R=R_{M',M}$ for each pair of modules $M',M$ in $\Cal O$.  
\b \m
\head \S 6.  Isomorphisms among quantum groups \endhead \m We will now 
investigate isomorphisms among the two-parameter quantum groups, and 
their connections with multiparameter quantum groups.  The case $n=2$ 
is special, and the two parameters collapse to one in the following 
sense.  Let $r,r',s,s'\in \K^\#$ and $r\neq s, \ r'\neq s'$.  If 
$rs^{-1}= r'(s')^{-1}$, there is an isomorphism of Hopf algebras
$$\phi: U_{r,s}({\frak {sl}}_2)\rightarrow U_{r',s'}({\frak {sl}}_2)$$
given by $\phi(\omega^{\pm 1})=\vo^{\pm 1}, \ \phi((\omega')^{\pm 1})
=(\vo ')^{\pm 1}, \ \phi(e)=\ve, \ \phi(f)= r^{-1}r'\vf$, where
``$\ \check{}\ $'' denotes generators of $U_{r',s'}({\frak {sl}}_2)$.  (When dealing with $\frak{sl}_2$, we 
omit the subscript ``1'' on the generators.) The proof is a simple 
check that the relations and coproducts are preserved.  In particular, 
if $q$ is a square root of $rs^{-1}$, then $U_{r,s}({\frak 
{sl}}_2)\cong U_{q,q^{-1}}({\frak {sl}}_2)$ as Hopf algebras.  
Therefore the one-parameter quantum group $U_q({\frak {sl}}_2)$ is 
isomorphic to the quotient of $U_{r,s}({\frak {sl}}_2)$ by the ideal 
generated by $\omega '-\omega^{-1}$.

If $n\geq 3$, there is no such isomorphism, as the following proposition
shows.  
\b
\proclaim{Proposition 6.1} Let $n\geq 3$, and assume there is an isomorphism
of Hopf algebras
$$\phi: U_{r,s}({\frak {sl}}_n){\rightarrow} U_{q,q^{-1}}
({\frak {sl}}_n)$$
for some $q$.  Then $r=q$ and $s=q^{-1}$.
\endproclaim
\m
\demo{Proof} Let $\phi$ be an isomorphism as hypothesized, and assume
$$\pi: U_{q,q^{-1}}({\frak {sl}}_n)\rightarrow U_q({\frak {sl}}_n)$$
is the surjection onto the standard one-parameter quantum group of [Ja]
given by $\pi(e_i)=E_i$, $\pi(f_i)=F_i$, $\pi(\omega_i^{\pm 1})=
K_i^{\pm 1}$, $\pi((\omega_i')^{\pm 1})=K_i^{\mp 1}$. 
For $1\leq i\leq n-1$, we have
$$\aligned
\Delta (\pi\phi (e_i)) = \pi\phi(\Delta(e_i)) &=\pi\phi(e_i\otimes 1 +
\omega_i\otimes e_i)\\
  &= \pi\phi(e_i)\otimes 1 + \pi\phi(\omega_i)\otimes \pi\phi(e_i).
\endaligned \tag 6.2 $$
Note that
$\pi\phi(\omega_i)$ is necessarily a group-like element.
Therefore $\pi\phi(e_i)$ is a skew-primitive element in $U_q({\frak {sl}}_n)$.
By Theorem 5.4.1, Lemma 5.5.5, and the subsequent comments in [M], the set
of group-like elements in $U_q({\frak {sl}}_n)$ is the group $G$
generated by $K_j$ $(1\leq j\leq n-1)$, and the 
skew-primitive elements together with the 
group-like elements span the subspace
$$\sum_{j=1}^{n-1} (\K E_j + \K F_j) +\K G.$$ 
Therefore 
$$\pi\phi(e_i)= \sum_{j=1}^{n-1} \alpha_j^i
E_j +\beta^i_j F_j + \sum_{g\in G}\gamma^i_g g$$
for some scalars $\alpha^i_j, \beta^i_j, \gamma^i_g \ (1\leq i,j\leq n-1, \
g\in G)$.  Consequently
$$\Delta(\pi\phi(e_i))=\sum_{j=1}^{n-1}
\alpha^i_j(E_j\otimes 1 + K_j\otimes E_j) + \beta_j^i (1\otimes F_j
+F_j\otimes K_j^{-1}) + \sum_{g\in G}\gamma^i_g g\otimes g.$$
By (6.2), this must be equal to 
$$\sum_{j=1}^{n-1} (\alpha^i_j E_j\otimes 1 +\beta^i_j F_j\otimes 1 + 
\pi\phi(\omega_i)\otimes \alpha^i_j E_j
+\pi\phi(\omega_i)\otimes \beta^i_j F_j) +
\sum_g(\gamma^i_g g\otimes 1 + 
\pi\phi(\omega_i)\otimes\gamma^i_g g).$$
Comparing these two expressions, and noting that $\pi\phi(\omega_i)\neq 1$
as $\phi$ is an isomorphism (or by the comparison of expressions), we see
first that
$$\sum_{j=1}^{n-1}\alpha^i_j E_j + \gamma^i_1 1 =
\sum_{j=1}^{n-1}( \alpha^i_j E_j + \beta^i_j F_j) +\sum_{g\in G}
\gamma^i_g g + \gamma^i_1 \pi\phi(\omega_i).$$
Therefore $\gamma^i_g=0$ for all $g$ except $g\in \{1, \pi\phi(\omega_i)\}$, \
$\gamma^i_{\pi\phi(\omega_i)}=-\gamma_1^i$, and all $\beta^i_j=0$.  We 
now have $$\pi\phi(e_i)=\sum_{1\leq j\leq n-1}\alpha^i_j E_j + 
\gamma^i_1
(1-\pi\phi(\omega_i)).$$
A further comparison of coproducts yields
$$\aligned \sum_{j=1}^{n-1} & \alpha^i_j \Big (E_j\otimes 1  + 
K_j\otimes E_j\Big) +\gamma^i_1\Big(1\otimes 1 -\pi\phi(\omega_i)\otimes 
\pi\phi(\omega_i)\Big) \\
& =\sum_{j=1}^{n-1}\alpha^i_j \Big( E_j\otimes 1 + 
\pi\phi(\omega_i) \otimes E_j\Big) \\
& \qquad + \gamma^i_1\Big(1\otimes 1 -\pi\phi(\omega_i)\otimes 1\Big) + 
\gamma^i_1\Big(\pi\phi(\omega_i)\otimes 1 -\pi\phi(\omega_i)\otimes
\pi\phi(\omega_i)\Big), \endaligned $$
which implies that
$$\alpha^i_j(K_j - \pi\phi(\omega_i))=0, \quad \text{for all} \ \  1\leq j 
\leq n-1.$$
Thus all $\alpha^i_j=0$ except possibly one, and some $\alpha^i_{j_i}$ 
must be nonzero as $\pi\phi$ is surjective.  Therefore 
$\pi\phi(\omega_i) =K_{j_i}$.

Next we will apply $\pi\phi$ to relation (R2') and use the relations in
$U_q({\frak {sl}}_n)$:
$$\aligned
\pi\phi(\omega_ie_i) &= \pi\phi(rs^{-1}e_i\omega_i)\\
K_{j_i}(\alpha^i_{j_i}E_{j_i}+ \gamma^i_1(1-K_{j_i})) & = rs^{-1}
(\alpha^i_{j_i} E_{j_i} + \gamma^i_1(1-K_{j_i}))K_{j_i}\\
\alpha^i_{j_i}q^2 E_{j_i}K_{j_i} + \gamma^i_1K_{j_i} - \gamma^i_1 K_{j_i}^2 &=
\alpha^i_{j_i} rs^{-1}E_{j_i}K_{j_i} + \gamma^i_1 rs^{-1}K_{j_i} - \gamma^i_1
rs^{-1} K_{j_i}^2.
\endaligned $$
This forces $rs^{-1} =q^2$, and $\gamma_1^i=0$ as $r\neq s$.  Applying
relation (R2') again, we have
$$\aligned
\pi\phi(\omega_ie_{i+1}) & = \pi\phi(se_{i+1}\omega_i)\\
K_{j_i} \alpha^{i+1}_{j_{i+1}} E_{j_{i+1}} & = \alpha^{i+1}_{j_{i+1}} s
E_{j_{i+1}} K_{j_i}\\
\alpha^{i+1}_{j_{i+1}} q^{\langle \alpha_{j_i},\alpha_{j_{i+1}}\rangle }
E_{j_{i+1}}K_{j_i} & = \alpha^{i+1}_{j_{i+1}} sE_{j_{i+1}}K_{j_i}.
\endaligned $$
Thus $s=q^{\langle \alpha_{j_i},\alpha_{j_{i+1}}\rangle }$ must hold.  By 
relations (R5) and (R6), we see that $| j_i-j_{i+1}|=1$, so that in 
fact $s=q^{-1}$.  Combining this with $rs^{-1}=q^2$ we now have $r=q$ 
as well.  \qed \enddemo \m \subhead Multiparameter quantum groups 
\endsubhead \m \m Deformations of $GL_n$ involving $1+\binom{n}{2}$ 
parameters were constructed independently by several authors (see 
[AST], [R], [S]).  The dual version is a multiparameter universal 
enveloping algebra, which was studied by Chin and Musson [ChM] and 
Dobrev and Parashar [DP].  We will show that our two-parameter quantum 
groups are essentially special cases of these multiparameter quantum 
groups, as should be expected.  We adopt the notation of Chin and 
Musson.

The $1+\binom{n}{2}$ parameters are denoted $\l$ and $p_{i,j} \ (1\leq i
<j\leq n)$ in [ChM].  Set $\l=rs^{-1}$ and $p_{i,j}=s^{-1}$ for all $i<j$.  
Let $\widehat{U}$ be the Hopf algebra generated by $E_i,F_i \ (1\leq 
i<n)$ and $K_i^{\pm 1},L_i^{\pm 1} \ (1\leq i \leq n)$ with relations 
given by \m 
\n {(ChM0)} \ \ The $K_i^{\pm 1}, L_j^{\pm 1}$ all 
commute with one another and $K_iK_i^{-1}=L_iL_i^{-1}$.  \m
 \n 
{(ChM1)} \ \ $K_jE_i=r^{-\delta_{i,j}}s^{-\delta_{i,j-1}}E_iK_j$ and 
$K_jF_i=r^{\delta_{i,j}}s^{\delta_{i,j-1}}F_iK_j$.  \m \n {(ChM2)} \ \ 
$L_jE_i=r^{\delta_{i,j-1}}s^{\delta_{i,j}} E_iL_j$ and 
$L_jF_i=r^{-\delta_{i,j-1}}s^{-\delta_{i,j}}F_iL_j$.  \m \n {(ChM3)} \ 
\ $E_iF_i -r^{-1}sF_iE_i =(r^{-1}s -1)(L_{i+1}K_{i+1} L_i^{-1}K_i^{-1} 
-1)$.  \m \n {(ChM4)} \ \ 
$E_iF_j=r^{\delta_{i,j+1}}s^{-\delta_{i,j-1}}F_jE_i$ if $i\neq j$.  \m 
\n {(ChM5)} \ \ $\ad_{\ell}(E_i)^{1-\langle 
\alpha_i,\alpha_j\rangle}(E_j)=0$ and 
$\ad_{\ell}(F_i)^{1-\langle\alpha_i,\alpha_j\rangle }(F_j)=0$ if 
$i\neq j$.  \b These relations are given in [ChM, Thm.  4.8] (for more 
general $\l, \ p_{i,j}$) as relations for a Hopf algebra that is 
defined as a subalgebra of the finite dual $A^0$ of a multiparameter 
quantum function algebra $A$.  In addition, Chin and Musson give one 
more set of conditions: (ChM6) \ those relations among the $K_i,L_j$ 
which determine the structure of the group they generate {\it as a 
subgroup of the group of units of} $A^0$.  This results in a Hopf 
algebra $\overline{U}$ (denoted $U$ in their paper).  Thus the 
multiparameter Hopf algebra $\overline{U}$ of Chin and Musson is the 
quotient of our $\widehat{U}$ by their relations (ChM6) (in case 
$\l=rs^{-1}$ and $p_{i,j}=s^{-1}$).

The Hopf structure of $\widehat{U}$ is defined by requiring $K_i, L_i$ to 
be group-like elements and $$\Delta(E_i)=E_i\otimes 1 + 
L_{i+1}L_i^{-1}\otimes E_i , \ \ \
\Delta(F_i)=F_i\otimes 1 + K_{i+1}K_i^{-1}\otimes F_i.$$
\m
\proclaim{Proposition 6.3}  There is a Hopf algebra morphism 
$\phi : \widehat{U} \rightarrow U_{r,s}({\frak {gl}}_n)$ given by
$$\aligned
\phi(L_i) & = a_1\cdots a_{i-1}b_{i+1}^{-1}\cdots b_n^{-1},\\
\phi(K_i) & = b_1^{-1}\cdots b_{i-1}^{-1}a_{i+1}\cdots a_n,\\
\phi(E_i) & = -s^{-1}(r-s)^2 e_i,\\
\phi(F_i) & = (\omega_i')^{-1}f_i.
\endaligned $$
\endproclaim
\m
The proof is just a check that the relations of $\widehat{U}$ are preserved
by $\phi$, and that the coproducts of the generators correspond to
the coproducts of their images.
\b
\n{\bf Remark 6.4.}
Note that the kernel and
cokernel of $\phi$ are generated by group-like elements.
We have $\omega_i=\phi(L_i^{-1}L_{i+1})$ and 
$\omega_i'=\phi(K_iK_{j+1}^{-1})$, which implies that $U_{r,s}({\frak {sl}}_n)$
is contained in the image of $\phi$.
In case $n=2$, it is easy to see that $\phi$ is an isomorphism.  
If $n=3$, straightforward calculations show that the cokernel of $\phi$ 
is precisely $\K\langle a_1^{\pm 1}\rangle$.
 
\b
\m

\vskip 3 mm
 
\address 
\newline 
Department of Mathematics, University of Wisconsin, Madison,
Wisconsin 53706   \newline 
benkart\@math.wisc.edu
\newline 
 \newline  
 Department of Mathematics, 
University of Massachusetts, Amherst, Massachusetts 01003  
\newline (2001-02) Department of Mathematics and Computer Science,
Amherst College, Amherst, Massachusetts 01002
\newline wither\@math.umass.edu \endaddress

\enddocument
\end